\documentclass[11pt,leqno]{article}

\usepackage{amsfonts,latexsym, amsmath, amssymb, amsthm, hyperref}

\usepackage{gfsartemisia}
\usepackage[T1]{fontenc}

\usepackage{color}
\usepackage{fancyhdr} 
\usepackage{lastpage}

\newtheorem{theorem}{Theorem}[section]
\newtheorem{lemma}[theorem]{Lemma}
\newtheorem{fact}[theorem]{Fact}

\newtheorem{proposition}[theorem]{Proposition}
\newtheorem{corollary}[theorem]{Corollary}

\theoremstyle{definition}

\newtheorem{remark}[theorem]{Remark}
\newtheorem{remarks}[theorem]{Remarks}
\newtheorem{quest}[theorem]{Question}
\newtheorem{note}[theorem]{Note}

\newtheorem*{note*}{Note}

\makeatletter

\numberwithin{equation}{section}

\makeatother

\newcommand{\prend}{$\hfill \quad \Box$}

\newcommand\blfootnote[1]{%
  \begingroup
  \renewcommand\thefootnote{}\footnote{#1}%
  \addtocounter{footnote}{-1}%
  \endgroup
}

\newcommand{\vertiii}[1]{{\left\vert\kern-0.25ex\left\vert\kern-0.25ex\left\vert #1 
    \right\vert\kern-0.25ex\right\vert\kern-0.25ex\right\vert}}

\begin{document}

\small

\title{Dichotomies, structure, and concentration in normed spaces}

\author{Grigoris Paouris\thanks{Supported by the NSF
CAREER-1151711 grant and Simons Foundation (grant \#527498).} \, and \, Petros Valettas\thanks{Supported by the NSF grant DMS-1612936.}}


\maketitle

\begin{abstract}\footnotesize
We use probabilistic, topological and combinatorial methods to establish the following deviation inequality:
For any normed space $X=(\mathbb R^n ,\|\cdot\| )$ there exists an invertible linear map $T:\mathbb R^n \to \mathbb R^n$ with
\[ \mathbb P\left( \big| \|TG\| -\mathbb E\|TG\| \big| > \varepsilon \mathbb E\|TG\| \right) 
\leq C\exp \left( -c\max\{ \varepsilon^2, \varepsilon \} \log n \right),\quad \varepsilon>0, \]
where $G$ is the standard $n$-dimensional Gaussian vector and $C,c>0$ are universal constants. It follows 
that for every $\varepsilon\in (0,1)$ and for every normed space $X=(\mathbb R^n,\|\cdot\|)$
there exists a $k$-dimensional subspace of $X$ which is
$(1+\varepsilon)$-Euclidean and $k\geq c\varepsilon \log n/\log\frac{1}{\varepsilon}$. 
This improves by a logarithmic on $\varepsilon$ term the best previously known result due to G. Schechtman.
\end{abstract}


\section{Introduction}

The concentration inequality in Gauss' space states that for any Lipschitz map $f:\mathbb R^n\to \mathbb R$ 
with $|f(x)-f(y)|\leq L\|x-y\|_2$ for all $x,y\in \mathbb R^n$ one has
\begin{align} \label{eq:con-Gauss}
\mathbb P \left( |f(G)- \mathbb Ef(G) | > t \right) \leq 2\exp(-\tfrac{1}{2}t^2/L^2), \quad t>0,
\end{align} where $G$ is the standard $n$-dimensional Gaussian vector (for a proof the reader is referred to \cite{Pis-pmBS}; see
\cite{Mau} for the precise constants). 
This inequality is the prototype of what is called nowadays the concentration of measure phenomenon, one
of the most important ideas in modern probability theory. This fundamental tool was put forward in the local theory 
of normed spaces in early 70's by V. Milman. Applying \eqref{eq:con-Gauss} for 
a norm $\|\cdot\|$ on $\mathbb R^n$ we get
\begin{align}  \label{eq:con-norms}
\mathbb P \left( \big| \|G\|- \mathbb E\|G\| \big|>t \mathbb E\|G\| \right) \leq 2 \exp( - \tfrac{1}{2} t^2 k), \quad t>0,
\end{align} where $k=k(X)=k(B_X):= (\mathbb E\|G\| / b)^2$ is referred to as the {\it critical dimension} (or {\it Dvoretzky number}) 
of the normed space 
$X=(\mathbb R^n, \|\cdot\|)$ and $b=b(X)=b(B_X)$ is the Lipschitz constant of the norm $\|\cdot\|$, i.e. $b=\max\{ \|\theta\|: \|\theta\|_2=1\}$.
It is well known that the above estimate is sharp in the large deviation regime, namely
\begin{align} \label{eq:ldr-lb}
\mathbb P( \|G\| \geq (1+t) \mathbb E\|G\|) \geq c \exp(-Ct^2k), \quad t\geq 1,
\end{align} where $c,C>0$ are universal constants\footnote{Here and everywhere else $C,c,C_1,c_1,\ldots$ stand for positive 
universal constants whose values may change from line to line. For any two quantities $A,B$ depending on the dimension, 
on the parameters of the problem, etc. 
we write $A\simeq B$ if there exists a universal constant $C>0$ -independent of everything- such that $A\leq C B$ and $B\leq C A$.} 
(see e.g. \cite[Corollary 3.2]{LedTal-book}, \cite[Statement 3.1]{LMS} and \cite[Proposition 2.10]{PVZ}). In 
the small deviation regime $0<t<1$ there exist many important examples which show that the obtained bounds are suboptimal;
see \cite{PVZ} and \cite{V-tight} for a detailed discussion. 
Ideally one would like to know what properties of the underlying function improve the concentration. An example of such a
result was recently obtained by the authors in \cite{PV-sdi} where they proved that a one-sided, variance-sensitive 
Gaussian small deviation inequality is valid for all convex functions.

This work is also concerned with optimal forms of the Gaussian 
concentration but the main focus will be on norms. Before stating the main problem of study, let us try to motivate the question
which describes it. It is known (see e.g. \cite{Mil-dvo}, \cite{Pis-pmBS}) that for any norm $\|\cdot\|$ on $\mathbb R^n$, 
there exists a $T\in GL(n)$ such that
\begin{align} \label{eq:con-John}
\mathbb P \left( \big| \| T G\|- \mathbb E\| TG\| \big|>t \mathbb E\| TG\| \right) \leq C \exp( -c t^2 \log n), \quad t>0,
\end{align} where $G\sim N({\bf 0},I_n)$. This follows from the fact that there exists a {\it position} (i.e. an invertible linear image) 
$T^{-1}(B_X)$ of the unit ball
$B_X= \{x\in \mathbb R^n: \|x\|\leq 1\}$ for which the critical dimension $k(T^{-1}B_X)$ is at least as large as $\log n$ (note that by
definition one has $c\leq k(X)\leq n$, hence a linear transformation is required to avoid degeneracy) and in turn this is 
combined with the general estimate
\eqref{eq:con-norms}. Traditionally, this 
is achieved for John's position, that is the maximal volume ellipsoid inscribed in $B_X$ is the standard Euclidean ball 
$B_2^n=\{ x\in \mathbb R^n: \|x\|_2\leq 1\}$. 
Let us note that the $\log n$ factor is optimal 
since for the $\ell_\infty$ norm 
we have that the cube $B_\infty^n=\{x\in \mathbb R^n : \|x\|_\infty \leq 1\}$ is in John's position and $k(\ell_\infty^n) \simeq \log n$. On the other hand,
the $\ell_\infty$ norm exhibits exponential concentration (see \cite{Tal-new-iso} and \cite{Sch-cube}),
\begin{align}
ce^{-C\varepsilon \log n} \leq \mathbb P \left( \big| \|G\|_\infty - \mathbb E\|G\|_\infty \big|>\varepsilon \mathbb E\|G\|_\infty \right) \leq C e^{- c\varepsilon \log n}, \quad 0<\varepsilon <1.
\end{align} In view of the above remarks the following question arises naturally:

\begin{quest} \label{que:main}
Is it true that for every norm $\|\cdot\|$ on $\mathbb R^n$ there exists a $T\in GL(n)$ with
\begin{align*}
\mathbb P \left( \big| \| TG\|- \mathbb E\| TG\| \big|>t \mathbb E\| TG\| \right) \leq C \exp( -c \max\{t^2,t\} \log n), \quad t>0 ?
\end{align*}
\end{quest}

One of the main difficulties to establish good small deviation estimates, is that the problem is ``isometric" in nature rather 
than ``isomorphic".

It is easy to see that two equivalent norms have Dvoretzky's number of the same order and therefore, by \eqref{eq:ldr-lb}, they exhibit the same large deviation
estimate. However, one may find two norms on $\mathbb R^n$ which are $2$-equivalent and the variance of one is polynomially small while the variance
of the other is only logarithmically small (with respect to the dimension), see for example \cite{PVZ}, \cite{PV-var} and \cite{LT}. 
We should mention that when the norm under consideration is close to the $\ell_2$-norm then it automatically exhibits the optimal 
concentration in terms of $\varepsilon$ and $n$ (Section \ref{sec:rems}, \S \, 1). On the other hand, 
for norms close to the $\ell_\infty$ norm this is no longer true: there exists a norm on $\mathbb R^n$ which is 2-equivalent to $\|\cdot\|_\infty$ and the optimal
concentration is as in \eqref{eq:con-norms} or \eqref{eq:con-John}, see e.g. \cite{PV-var}. In fact for any norm, one can construct a 5-equivalent norm 
for which \eqref{eq:con-norms} is sharp, see \cite[Section 3]{V-tight}. After these observations it seems quite discouraging to tackle Question \ref{que:main} and leaves only the hope that the appropriate selection of the linear transformation will fix the problem.

The above discussion leads naturally to the question of how we successfully select the position to establish improved concentration 
estimates. It turns out that John's position, which was widely used to attack problems lying in concentration estimates 
in the local theory of normed spaces, is not the proper one. It is possible to construct $n$-dimensional normed spaces whose unit 
ball is in John's position, the critical dimension can be of any order in the admissible range for this position and yet the optimal 
concentration is as in \eqref{eq:con-norms} (see Section \ref{sec:rems} for the exact formulation). In the closely related setting 
of almost Euclidean subspaces, it was emphasized by the authors in \cite{PV-dvo-Lp} that the choice of the position is crucial for
improving the estimates and subsequently, it was established by K. Tikhomirov in \cite{Tik-unc} that John's position may give suboptimal bounds.

To the positive direction, there are recent results which indicate that better concentration may be obtained for special 
classes of normed spaces by selecting a different position than John's; see \cite{PV-dvo-Lp} and \cite{Tik-unc}. 
In particular, in \cite{PV-dvo-Lp} for subspaces of $L_p, \; 1\leq p<\infty$ the authors showed that in Lewis' position 
(see \cite{TJ-book} for the related definition) the concentration is at least as good as for the $\ell_p$ norms and in
view of \cite{PVZ} this is best possible. Thus, it follows that for every $n$, for every $1\leq p\leq c\log n$ and for
every $n$-dimensional subspace $X$ of $L_p$ the answer in Question \ref{que:main} is affirmative. In \cite{Tik-unc} K. Tikhomirov
proved that for 1-unconditional norms in $\ell$-position (see \cite{TJ-book} for the definition) the concentration is at least as good 
as for the $\ell_\infty$ norm, thus establishing Question \ref{que:main} in affirmative for those spaces.
Let us mention that in all the above approaches, Gaussian functional 
inequalities are used rather than the classical Gaussian concentration in terms of the Lipschitz constant.

The first main step to tackle Question \ref{que:main}, is to show that the latter has an affirmative
answer when the normed space $X=(\mathbb R^n,\|\cdot\|)$ is not {\it extremal} with respect to the unconditional structure. In order to formulate the 
result, we have to recall some standard terminology. 
Following the notation used in Figiel and Johnson \cite{FJ}, for any basis $\{b_i\}_{i=1}^n$ in a normed space $X=(\mathbb R^n, \|\cdot\|)$
the {\it unconditional basis constant} is given by
\begin{align*}
{\rm unc}\{b_i\}_{i=1}^n =\sup \left\{ \left\| \sum_{i=1}^n \varepsilon_i t_ib_i \right\| :  \left\| \sum_{i=1}^n t_ib_i\right\| =1\right\},
\end{align*} where the supremum is over all choices of signs $\varepsilon_i=\pm 1$ and all scalars $(t_i)$. Next, we define
\begin{align*}
{\rm unc}X= \inf \left\{ {\rm unc}\{b_i\}_{i=1}^n : \{b_i\}_{i=1}^n \; \textrm{ is a basis for} \; X \right \}.
\end{align*} For any normed space $X=(\mathbb R^n,  \|\cdot\|)$ one has ${\rm unc}X\leq \sqrt{n}$. Indeed; note that
${\rm unc}X$ is nothing more than the Banach-Mazur distance of $X$ to the class of all $n$-dimensional $1$-unconditional normed spaces and 
hence, ${\rm unc}X \leq d(X,\ell_2^n)\leq \sqrt{n}$, by John's theorem \cite{Jo}. At this point we should mention that this bound is optimal
(up to universal constants), since there exists $n$-dimensional normed space $E$ 
(in fact any ``typical" subspace of $\ell_\infty^n$ of proportional dimension) which satisfies 
${\rm unc}E\geq c \sqrt{n}$ (see e.g. \cite{FKP} and \cite{FJ}). With this notation we have the following:
 
 \begin{theorem} \label{thm:main-1}
 Let $X=(\mathbb R^n,\|\cdot\|)$ be a normed space. Then, there exists a $T\in GL(n)$ such that
 \begin{align*}
 \mathbb P \left( \big| \| TG\|- \mathbb E\| TG\| \big|>t \mathbb E\| TG\| \right) \leq 
 C \exp \left( -c \max\{t^2,t\} \log \left( \frac{e n}{ ({\rm unc}X)^2 } \right) \right), \quad t>0 ,
 \end{align*} where $G$ is the standard $n$-dimensional Gaussian vector.
 \end{theorem} The proof of Theorem \ref{thm:main-1} uses probabilistic and topological tools. For the probabilistic part we use 
 Talagrand's improvement upon the classical Gaussian Poincar\'e inequality.
 
\begin{theorem} [Talagrand's $L_1-L_2$ bound] \label{thm:Tal-L1-L2} 
For any absolutely continuous function $f:\mathbb R^n \to \mathbb R$, we have
\begin{align} \label{eq:Tal-L1-L2}
{\rm Var}_{\gamma_n} (f) \leq C \sum_{i=1}^n \frac{\|\partial_i f\|_{L_2}^2}{ 1+ \log \left(\|\partial_i f\|_{L_2} / \|\partial_i f\|_{L_1} \right)},
\end{align} where $\partial_i f$ stands for the i-th partial derivative of $f$.
\end{theorem}

Talagrand proved the above theorem for the Hamming cube \cite{Tal-russo} and the Gaussian version of this inequality was 
presented in \cite{CorLed}. It is known that this inequality is also related to the {\it super-concentration} phenomenon. 
Following Chatterjee, a function $f:\mathbb R^n\to \mathbb R$ is said to be $\varepsilon_n$-super-concentrated if 
\[ {\rm Var}[f(G)] \leq \varepsilon_n \mathbb E \|\nabla f(G)\|_2^2. \]
The reader is referred to \cite{Cha} for a detailed exposition of this very interesting subject and further applications. 
The use of Talagrand's inequality in the study of the 
asymptotic theory of finite-dimensional normed spaces was put forward by the authors and J. Zinn in \cite{PVZ}. 
The authors there, use the inequality to prove sharp concentration for the $\ell_p$ norms when $p$ grows along with $n$. 
Additionally, it is proved that the $\ell_p$ norms are super-concentrated for $p> c \log n$.
This inequality also played a central role in the aforementioned work of K. Tikhomirov \cite{Tik-unc} for the case of 1-unconditional norms.

As long as the choice of the position is concerned, we should mention that most canonical positions that are in use in the geometry of 
finite-dimensional normed spaces arise as solution of an extremum value problem, where a geometric functional is optimized subject
to a constraint. For a detailed exposition of this positions the reader is referred to \cite{GM-iso}. 
However, in our approach, we prove the existence of a position in Theorem \ref{thm:main-1} by employing 
a topological tool, namely the Borsuk-Ulam antipodal theorem \cite{Bor}, \cite{Mat}. The Borsuk-Ulam theorem has already found many
fruitful applications in geometric (linear and nonlinear) functional analysis, see e.g. \cite{L}, \cite{R}, \cite{GM}, \cite{KL} just to name a few. 

The unconditional constant in Theorem \ref{thm:main-1} is naturally involved, since inequality \eqref{eq:Tal-L1-L2} is sensitive
with respect to the coordinate structure. One can find good local unconditional structure based on a fundamental result of Alon 
and Milman \cite{AM}. The latter states that for every $\varepsilon\in (0,1)$ 
there exists a constant $C(\varepsilon) >0$ such that every normed space 
$X=(\mathbb R^n, \|\cdot\|)$ satisfies the following dichotomy:
\begin{itemize}
\item Either there exists a subspace $E$ with $\dim E=m\geq e^{\sqrt{\log n}}$ and $d(E, \ell_2^m)<1+\varepsilon$, 
\item Or there exists a subspace $F$ with $\dim F=k\geq e^{\sqrt{\log n}- e^{C(\varepsilon)}}$ and $d(F,\ell_\infty^k)< 1+\varepsilon$.
\end{itemize}

The key tool, the authors prove in \cite{AM}, for establishing the above dichotomy is a combinatorial result for locating 
$\ell_\infty$-structure (see Section \ref{sec:dicho} for the precise formulation). Combining this with Theorem 
\ref{thm:main-1} we prove the following probabilistic dichotomy for the Gaussian concentration.

\begin{theorem} \label{thm:main-2}
Let $X=(\mathbb R^n, \|\cdot\|)$ be a normed space whose unit ball $B_X$ is in John's position and let $0< \delta <1/2$. 
Then, we have the following dichotomy:
\begin{itemize}
\item Either the random\footnote{Here the randomness is considered with respect to the unique probability measure on the 
Grassmannian $G_{n,k}$ which is invariant under the orthogonal group action.} subspace $E$ with $\dim E=k \geq n^{1/2-\delta}$ satisfies
\begin{align*}
\mathbb P\left( \big| \|G\|_{E\cap B_X} -\mathbb E\|G\|_{E\cap B_X} \big| > t \mathbb E\|G\|_{E\cap B_X} \right) \leq C e^{-c t^2 k}, \quad t>0,
\end{align*}
\item Or there exists a subspace $F$ with $\dim F=m\geq c n^{1/2}$ and an invertible linear map $T:F \to F$ such that
\begin{align*}
\mathbb P \left( \big| \|TG\|_{F\cap B_X} - \mathbb E\|TG\|_{F\cap B_X} \big| >t\mathbb E\|TG\|_{F\cap B_X} \right) \leq C e^{-c\delta \max \{ t^2,t \} \log m }, \quad t>0,
\end{align*} 
\end{itemize} where $G$ is the standard Gaussian vector and $c,C>0$ are universal constants.
\end{theorem}

The above result can be interpreted as the probabilistic aspect of Alon-Milman theorem for the Gaussian measure. 
In turn, this is sufficient to imply an affirmative answer to Question \ref{que:main}.

\begin{theorem}  \label{thm:main-3}
Let $X=(\mathbb R^n, \|\cdot\|)$ be a normed space. There exists a $T\in GL(n)$
such that for all $t>0$,
\begin{align*}
\mathbb P\left( \big| \|TG\| -\mathbb E\|TG\| \big| > t \mathbb E\|TG\| \right) \leq C \exp ( -c\max\{t^2,t\} \log n ), 
\end{align*} where $G$ is the standard Gaussian vector.
\end{theorem} 
We would like to emphasize the fact that the position involved in the above theorem is not one of the standard positions used in the local theory of
Banach spaces. It is rather a position derived by the combination of all the aforementioned techniques. It would be interesting to know that Theorem 
\ref{thm:main-3} holds true for some classical position, such as the position of minimal $M$ or the $\ell$-position (see Section 2 for the related definitions; see also
Section 5 for the fact that this cannot be achieved in John's position). 

The Gaussian concentration for norms is closely related to the local almost Euclidean structure. In his seminal work \cite{Mil-dvo}, V. Milman
establishes a random version of the celebrated result of Dvoretzky \cite{Dvo} on the almost spherical sections of convex bodies; 
see also \cite{FLM}, \cite{MS}. V. Milman uses \eqref{eq:con-norms} to prove that for any $\varepsilon\in (0,1)$
the random $m$-dimensional subspace $E$ of $X$ (with respect to the Haar measure on the Grassmannian) is $(1+\varepsilon)$-spherical, i.e.
\begin{align*}
(1-\varepsilon) a B_E \subset B_X\cap E \subset (1+\varepsilon) a B_E
\end{align*} for some appropriate constant $a>0$ depending only on $X$, with probability greater than $1-e^{-c\varepsilon^2 k(X)}$, 
as long as $m\leq c \varepsilon^2 k(X)$ (see \cite{Go} and
\cite{Sch-eps} for the dependence $\varepsilon^2$). Thus, if we define $k_r(X,\varepsilon)$ to be the maximal $k$ for which 
the {\it random} $k$-dimensional subspace of $X$ is $(1+\varepsilon)$-spherical with probability at least 2/3 say, Milman's argument shows that 
$k_r(X,\varepsilon)\geq c\varepsilon^2 k(X)$. Theorem \ref{thm:main-3} then, combined with a standard net argument, implies the following.

\begin{corollary} \label{cor:main-4}
For every normed space $X=(\mathbb R^n, \|\cdot\|)$ there exists a position $B$ of $B_X$ such that for every $\varepsilon\in (0,1)$ one has
$k_r(B,\varepsilon) \geq c\varepsilon  \log n / \log(1/\varepsilon)$.
That is, the random $k$-dimensional section of $B$ with $k\leq c\varepsilon \log n/ \log(1/\varepsilon)$ is $(1+\varepsilon)$-spherical with
probability greater than $1-n^{-c\varepsilon}$.
\end{corollary}

It is a question of Grothendieck \cite[\S \, 7]{Gro} to determine the largest possible $k=k(n, \varepsilon)$ for which 
{\it every} $n$-dimensional normed space
$X$ admits a $k$-dimensional subspace which is $(1+\varepsilon)$-Euclidean. More precisely, for given $\varepsilon \in (0,1)$ and 
$X=(\mathbb R^n, \|\cdot\|)$ we denote by $k(X,\varepsilon)$ the largest $k$ so that {\it there exists} a $k$-dimensional subspace of $X$
which is $(1+\varepsilon)$-Euclidean. Then, for $0<\varepsilon<1$ we set
\[ k(n,\varepsilon) = \inf \{ k(X,\varepsilon) : \dim X =n\} . \] 
It is easy to show (see e.g. \cite{Sch-surv}) that 
\[ k(n,\varepsilon) \leq k(\ell_\infty^n,\varepsilon) \leq C \log n / \log (1/ \varepsilon) , \] 
while the fundamental fact that the function $k(n, \varepsilon)\to \infty$, for $\varepsilon = \varepsilon_n \to 0$ as $n\to \infty$ has 
first been established by Dvoretzky in \cite[Theorem 1]{Dvo}, who showed the quantitative estimate 
$k(n,\varepsilon)\geq c\varepsilon \sqrt{\log n }/ \log\log  n$. The aforementioned randomized version of Dvoretzky's theorem by 
V. Milman \cite{Mil-dvo} improved the bound to $k(n, \varepsilon)\geq c \varepsilon^2 \log n /\log\frac{1}{\varepsilon}$ 
and then Gordon in \cite{Go} showed that $k(n, \varepsilon) \geq c\varepsilon^2 \log n$ 
(see also \cite{Sch-eps} for an alternative proof of this estimate). Schechtman proved in \cite{Sch-eps-exist}
that one can always have $k(n, \varepsilon)\geq c\varepsilon \log n /( \log\frac{1}{\varepsilon})^2$.
Corollary \ref{cor:main-4} also gives the best known estimate for this question up-to-date, but the best possible lower 
estimate for the function $k(n,\varepsilon)$ remains a fundamental open problem.

Our approach shares common points with Schechtman's argument but is essentially
different. In both cases the Alon-Milman theorem is a crucial tool. Schechtman invokes an iteration
scheme based on James' distortion lemma \cite{Ja} to find further a subspace which is sufficiently close to $\ell_\infty$. This strategy 
is followed because $\ell_\infty$ admits finer dependence on $\varepsilon$ for the existential Dvoretzky; see \cite{Sch-surv} and \cite{Sch-eps-exist}. 
This procedure yields a redundant logarithmic term of $\varepsilon$ compared to Corollary \ref{cor:main-4}. To the contrary, we use Alon-Milman theorem 
to determine the linear map $T$ for which the norm exhibits at least as good concentration as in the $\ell_\infty$ case. Thus, we obtain
the dependence on $\varepsilon$ that holds true for the random version of Dvoretzky's theorem in the case of $\ell_\infty$, see e.g. \cite{Sch-cube}
and \cite{Tik-cube}. 

Concluding, we would like to point out that V. Milman in \cite{Mil-observ} had observed the connection of topological tools with the
problem of the dependence on $\varepsilon$ in Dvoretzky's theorem, yet it hadn't been exploited until now.

The rest of the paper is organized as follows: In Section \ref{sec:2-level} we use Talagrand's $L_1-L_2$ bound 
to establish a two-level Gaussian deviation inequality for Lipschitz functions, where the Lipschitz condition 
is considered in both $\ell_2$ and $\ell_\infty$ sense. In Section
\ref{sec:con-d-unc} we present the proof of Theorem \ref{thm:main-1}. In Section \ref{sec:dicho} we employ Theorem \ref{thm:main-1}
and the Alon-Milman theorem to obtain Theorem \ref{thm:main-3}. Finally, in Section \ref{sec:rems} we conclude with remarks and 
questions that arise from our work. For background material on the geometry of Banach spaces the reader may consult the 
monographs \cite{MS, TJ-book, Pis-book, AGM}.

\section{A two-level Gaussian deviation inequality} \label{sec:2-level}

It is well known that Poincar\'e inequalities imply exponential concentration for Lipschitz maps (see \cite{Led-book} and \cite{BLM}). 
Since Talagrand's $L_1-L_2$ inequality is an improved version of the classical Poincar\'e inequality one gets straightforward 
improvements on the corresponding exponential concentration. In order to illustrate that, let us examine what is the corresponding 
deviation estimate we obtain by Talagrand's inequality, if we employ the standard method of bounding the variance of 
the moment generating function of 
a Lipschitz map $f$. To this end, we introduce the following notation: For any Lipschitz map $f:\mathbb R^n\to \mathbb R$ let
\begin{align*}
b=b(f):=\inf\{t>0 : |f(x)-f(y)|\leq t\|x-y\|_2 , \; \forall \, x,y\in \mathbb R^n \}.
\end{align*} and similarly
\begin{align*}
a=a(f):=\inf\{t>0 : |f(x)-f(y)|\leq t \|x-y\|_\infty, \; \forall \, x,y\in \mathbb R^n\}.
\end{align*}
Note that in the light of $\|\cdot\|_\infty \leq \|\cdot\|_2 \leq \sqrt{n} \|\cdot\|_\infty$ one has
\begin{align}
b\leq a\leq b\sqrt{n}.
\end{align} If $\|\cdot\|$ is an arbitrary norm on $\mathbb R^n$ and 
$|f(x)-f(y)|\leq L \|x-y\|$ for all $x,y\in \mathbb R^n$, then 
\begin{align} 
\|\nabla f(x)\|_\ast \leq L,
\end{align} where $\|\cdot\|_\ast$ is the dual norm of $\|\cdot\|$, i.e.
\begin{align*}
\|y\|_\ast= \sup\{ \langle x,y\rangle : \|x\|\leq 1\}, \quad y\in \mathbb R^n,
\end{align*} and the gradient of $f$ is defined almost everywhere by Rademacher's theorem, see e.g. \cite{EG}.

Now we come to the aforementioned improvement of the exponential concentration via the $L_1-L_2$ bound. 
In order to simplify considerably the computations let us assume that $f$ has some symmetries, 
i.e. $f$ is permutation invariant\footnote{For any permutation $\pi: [n]\to [n]$ we define the permutation matrix
$P_\pi$ associated with $\pi$ as follows: $P_\pi( e_i) =e_{\pi(i)}$. Note that $P_\pi \circ P_\sigma= P_{\pi\sigma}$ and 
$P_\pi^{-1}= P_{\pi^{-1}} = P_\pi^\ast$ for all permutations $\pi, \sigma$. 
A function $f: \mathbb R^n\to \mathbb R$ is said to be 
permutation invariant if $f\circ P_\pi =f$ for any permutation $\pi$.}. In that case we have that $h=e^{\lambda f}, \lambda >0$ is also permutation invariant and
\begin{align*}
\partial_i h = \partial_i (h \circ P_\pi) =  \langle P_{\pi}^\ast \circ {\nabla h} \circ P_\pi ,e_i \rangle = (\partial_{\pi(i)} h) \circ P_\pi, \quad i=1,2,\ldots,n,
\end{align*} for any permutation $\pi$. It follows that
\begin{align*}
\lambda^p \mathbb E e^{ p \lambda f}|\partial_i f|^p  = \|\partial_i h\|^p_{L_p} = \|\partial_{\pi (i)} h\|_{L_p}^p,
\end{align*} for all $i\leq n$, for any permutation $\pi$ and for $p>0$, since $P_\pi$ is orthogonal. In particular, the $L_2$-norm 
of all partial derivatives of $h$ are equal, thus
\begin{align*}
\|\partial_i h\|_{L_2}^2 =\frac{1}{n} \sum_{i=1}^n \|\partial_i h\|_{L_2}^2 =\frac{\lambda ^2}{n} \mathbb E e^{2\lambda f} \|\nabla f\|_2^2 \leq \frac{\lambda^2 b(f)^2}{n} \mathbb Ee^{2\lambda f}, \quad i=1,\ldots, n.
\end{align*} Arguing similarly, we get
\begin{align*}
\|\partial_i h\|_{L_1} \leq \frac{\lambda a(f)}{n} \mathbb Ee^{\lambda f}, \quad i=1,\ldots, n.
\end{align*} Applying Theorem \ref{thm:Tal-L1-L2} for $h$ and taking into account the previous estimates we obtain
\begin{align*}
\mathbb Ee^{2\lambda f}- (\mathbb Ee^{\lambda f})^2= {\rm Var}(e^{\lambda f}) 
\leq  \frac{ C \lambda ^2 b^2}{1+ \log (nb^2/a^2)} \mathbb E e^{2\lambda f}, \quad \lambda>0,
\end{align*} where $a=a(f)$ and $b=b(f)$. Next, we argue as in \cite[p.70]{BLM}. Set $\rho^2= \frac{C b^2}{1+\log(nb^2/a^2)}$ and
apply the previous estimate for $F=f-\mathbb Ef$ and $\lambda=s/\rho$. Then, for $\psi(s)= \mathbb Ee^{2s F/\rho}$, 
we obtain the following conditions:
\begin{align*}
\lim_{s \to 0} \frac{\psi(s)-1}{s}=0, \quad (1-s^2)\psi(s) \leq (\psi(s/2))^2, \quad  0<s<1.
\end{align*} It is a calculus exercise to show that such a function satisfies $\psi(s) \leq (1-s^2)^{-2}$ for all $0<s<1$. In particular, $\psi(1/\sqrt{2})\leq 4$.
Using Markov's inequality and combining with the 
classical Gaussian concentration \eqref{eq:con-Gauss} we conclude the following.

\begin{proposition}\label{prop:conc-permut}
Let $f:\mathbb R^n\to \mathbb R$ be a permutation invariant function. If $|f(x)-f(y)|\leq b \|x-y\|_2$ and 
$|f(x)-f(y)|\leq a \|x-y\|_\infty$ for all $x,y\in \mathbb R^n$, then 
\begin{align} \label{eq:con-perm}
\mathbb P \left( \left| f(G)-\mathbb E[f(G)] \right| > t \right) \leq 4 \exp \left( -c \max\left\{  \frac{t^2}{b^2}, \frac{t}{b} \sqrt{ \log (nb^2/a^2)} \right\} \right),
\end{align} for all $t>0$, where $G$ is the standard $n$-dimensional Gaussian vector and $c,C>0$ are universal constants.
\end{proposition}

The purpose of this section is to show that a similar concentration inequality can be proved regardless the symmetries of $f$.
To this end, we will need the following consequence of Talagrand's $L_1-L_2$ inequality (for a proof see e.g. \cite[Theorem 5.4]{Cha}):

\begin{lemma} \label{lem:key}
Let $f:\mathbb R^n \to \mathbb R$ be an absolutely continuous function and let 
\begin{align*}
R(f) =\frac{ \mathbb E \|\nabla f(G)\|_2^2}{ \sum_{i=1}^n (\mathbb E |\partial_i f(G)|)^2}.
\end{align*} Then, we have
\begin{align} \label{eq:key} 
{\rm Var}[f(G)] \leq C \frac{\mathbb E \|\nabla f(G)\|_2^2}{1+\log R(f) },
\end{align} where $G$ is the standard Gaussian vector on $\mathbb R^n$.
\end{lemma}

\medskip

Now we are ready to prove the aforementioned two-level deviation inequality. This inequality is in the spirit of
Talagrand's two-level deviation inequality for the exponential distribution \cite{Tal-new-iso} (see also \cite{BL} for an alternative proof).

\begin{proposition} \label{prop:2-level}
Let $f:\mathbb R^n\to \mathbb R$ be a Lipschitz map with
\begin{align*}
 |f(x)-f(y)|\leq b \|x-y \|_2 , \quad |f(x)-f(y)|\leq a \|x-y \|_{\infty}, \quad x,y\in \mathbb R^n
\end{align*} and $\| \partial_i f \|_{L_1} \leq A$ for all $i\leq n$. Then, if we set $F=f-\mathbb Ef$, for 
all $\lambda >0$ we have
\begin{align} \label{eq:log-mgf}
{\rm Var}(e^{\lambda F}) \leq \frac{C\lambda^2 b^2}{\log (e+\frac{ b^2}{a A})} \mathbb E e^{2\lambda F}.
\end{align} 
Moreover, we obtain
\begin{align} \label{eq:exp-conc}
\mathbb P \left ( | f(G)- \mathbb E[f(G)] | >t \right) 
\leq 4\exp \left(-c\max \left \{ \frac{t^2}{b^2},  \frac{t}{b} \sqrt{ \log \left( e +\frac{b^2}{aA} \right) } \right \} \right), \quad t>0,
\end{align} where $C,c>0$ are universal constants.
\end{proposition}

\noindent {\it Proof.} We fix $\lambda>0$ and we apply Lemma \ref{lem:key} for the function $e^{\lambda F}$. Then, we obtain
\begin{align*}
{\rm Var}(e^{\lambda F}) \leq C\lambda^2 \frac{\mathbb E (\|\nabla f\|_2^2 e^{2\lambda F}) }{1+ \log\left( \frac{\mathbb E \|\nabla f\|_2^2 e^{2\lambda F}}{w} \right) },
\end{align*} where $w=\sum_{i=1}^n \left(\mathbb E |\partial_i f|e^{\lambda F} \right)^2$.
Note that $\mathbb E\|\nabla f\|_2^2 e^{2\lambda F} \leq b^2 \mathbb E e^{2\lambda F}$ and the function
$z\mapsto \frac{z}{1+\log(z/w)} $ is non-decreasing for $z\geq w$, hence we get
\begin{align*}
{\rm Var}(e^{\lambda F}) \leq  \frac{C\lambda^2 b^2 }{1+\log \left( w^{-1} b^2 \mathbb Ee^{2\lambda F}\right) } \mathbb Ee^{2\lambda F}.
\end{align*}
Finally, note that 
\begin{align*}
w= \sum_{i=1}^n \left( \mathbb E |\partial_i f| e^{\lambda F} \right)^2 
\leq  \sum_{i=1}^n \mathbb E e^{\lambda F} \mathbb E|\partial_i f|^2 e^{\lambda F}  \leq  A \mathbb Ee^{\lambda F} \mathbb E\|\nabla f\|_1 e^{\lambda F} \leq 
aA (\mathbb E e^{\lambda F})^2 ,
\end{align*} where we have used the Cauchy-Schwarz inequality, the bound on the $L_1$ norm of the partial derivatives, and the pointwise bound $\|\nabla f(x)\|_1\leq a$. 
This completes the proof of the first assertion. The concentration estimate \eqref{eq:exp-conc} can be proved in a standard fashion as before, by using \eqref{eq:log-mgf}.
See e.g. \cite[p.70]{BLM} for the details. \prend

\begin{remark} The logarithm appearing on the estimate \eqref{eq:exp-conc} is almost the same as in \eqref{eq:con-perm} without the parameter
$A$. Note that in general the least possible $A$ satisfies $A\leq a$. In the next section we show that after composing the 
function with a suitable diagonal matrix we may bound $A\leq a/n$ and hence derive exactly the estimate \eqref{eq:con-perm}.
\end{remark}

Next, we present an application of the previous distributional inequality in the context of $1$-unconditional norms. To this end
let us recall the definition of the {\it position of minimal $M$}: A norm $\|\cdot\|$ on $\mathbb R^n$ is said to be in position of minimal $M$ if
for every $T\in SL(n)$ we have
\begin{align*}
\mathbb E\|G\| \leq \mathbb E\|TG\|, \quad G\sim N({\bf 0},I_n).
\end{align*} In this case the norm satisfies the following {\it isotropic condition}:
\begin{align} \label{eq:iso}
\int_{\mathbb R^n} \langle \nabla \| x\|, \theta \rangle \langle x, \theta \rangle \, d\gamma_n(x) =\frac{ \mathbb E\|G\|}{n}, \quad \theta\in S^{n-1}.
\end{align} For more properties of this position the reader is referred to \cite{GM-iso}.

\begin{proposition} \label{prop:con-unc}
Let $\|\cdot\|$ be a 1-unconditional norm on $\mathbb R^n$ which is in position of minimal $M$. Then, we have the following 
distributional inequalities:

\begin{itemize}

\item [\rm (i)] For all $t>0$,
\begin{align} \label{eq:uncon-conc-1}
\mathbb P \left( \big| \|G\| - \mathbb E\|G\| \big| > t \mathbb E\|G\| \right) \leq C\exp \left(- c \max \left \{t^2 k, t \sqrt{k\log(en/k)} \right\} \right), 
\end{align} where $k=k(X)$ and $X=(\mathbb R^n, \|\cdot\|)$. 

\item [\rm (ii)] In particular $k\geq c\log n$, hence
\begin{align} \label{eq:uncon-conc-2}
\mathbb P \left( \big|  \|G\| - \mathbb E\|G\| \big| > t \mathbb E\|G\| \right) \leq C\exp \left(- c \max \{t^2, t \} \log n \right), \quad t>0,
\end{align} 
\end{itemize}
where $G$ is the standard $n$-dimensional Gaussian vector and $C,c>0$ are universal constants.
\end{proposition}

\noindent {\it Proof.} (i). Set $f(x)=\|x\|$. The unconditionality and the convexity of $f$ implies that $x_i \mapsto \partial_if(x) $ is nondecreasing
function of $|x_i|$, Hence, Chebyshev's association inequality (see e.g. \cite[Section 2.10]{BLM}) yields
\begin{align*}
\mathbb E |\partial_i f(G)| \cdot \mathbb E|g_i| \leq \mathbb E |g_i \partial_i f(G)| = \frac{\mathbb Ef(G)}{n} \Longrightarrow \|\partial_i f\|_{L_1(\gamma_n)} \leq \frac{c_1 \mathbb Ef(G)}{n},
\quad i=1,\ldots,n.
\end{align*} Furthermore we have the following:

\smallskip

\noindent {\it Claim.} Note that $a(f)=\max\{ \|x\| : \|x\|_\infty \leq 1\}$ and $\mathbb E \|G\| \geq c a(f)$.

\smallskip

\noindent {\it Proof of Claim.} Indeed; we may write
\begin{align*}
\mathbb E \left\| \sum_{i=1}^n g_i e_i \right\|=\mathbb E_\varepsilon \mathbb E \left\| \sum_{i=1}^n \varepsilon_i |g_i| e_i \right\| 
\geq \mathbb E_\varepsilon \left \| \sum_{i=1}^n \varepsilon_i \mathbb E|g_i| e_i\right \|= \sqrt{ \frac{2}{\pi}} \left\|\sum_{i=1}^n e_i \right\|,
\end{align*} by Jensen's inequality and the unconditionality of the norm. On the other hand we have
\begin{align*}
\|i :\ell_\infty^m \to X\| =\max_{\|x\|_\infty\leq 1}\|x\| = \max_{\varepsilon_i=\pm 1}\left\| \sum_{i=1}^n \varepsilon_i e_i \right\| = \left \| \sum_{i=1}^n e_i \right \|,
\end{align*} which proves the assertion. \prend

\smallskip

\noindent Thus, a straightforward application of Proposition \ref{prop:2-level} combined with the above estimates 
yields \eqref{eq:uncon-conc-1}. 

\smallskip

\noindent (ii). From \eqref{eq:uncon-conc-1} applied for $t\simeq 1$ and 
compared with \eqref{eq:ldr-lb} we get $k\geq c\log(en/k)$ which yields the desired estimate. \prend

\begin{note}
In \cite{Tik-unc} K. Tikhomirov proves \eqref{eq:uncon-conc-2} using the $\ell$-position. The
latter is a variant of the position we use here and satisfies an analogous isotropic condition, namely
\begin{align} \label{eq:iso-ell}
\int_{\mathbb R^n} \langle \nabla \| x\|, \theta \rangle \langle x, \theta \rangle \|x\|\, d\gamma_n(x) =\frac{ \mathbb E\|G\|^2}{n}, \quad \theta\in S^{n-1}. 
\end{align} Using \eqref{eq:iso-ell} instead, we can again show that $n \|\partial_i f\|_{L_1}\leq c \mathbb E f(G)$ for all $i\leq n$, 
where $f(x)=\|x\|$ is 1-unconditional. Note that, according to the argument of Proposition \ref{prop:con-unc}.(i), it suffices to obtain an 
upper bound of the form $\mathbb E |g_i \partial_i f(G)| \leq C \mathbb Ef(G) /n$ for each $i\leq n$. We will show that indeed this follows from \eqref{eq:iso-ell} and the unconditionality of $f$. To this end recall the known fact (see e.g. \cite[Proposition 16]{Tik-unc} for a proof) that $k(X)\geq c_1\log n$ for any normed space $X=(\mathbb R^n, \|\cdot\|)$ whose unit ball $B_X$ is in $\ell$-position. Using the small ball probability estimate for norms 
(see e.g. \cite{LO}, \cite{KV} and \cite{PV-sdi} for a refinement) we have
\begin{align*}
A:=\{ \|G\| \leq c_0 \mathbb E\|G\| \} , \quad \mathbb P(A) \leq e^{-c_2 k(X)} \leq 1/n^2.
\end{align*} Thus, using the unconditionality, the fact that 
\[ |\partial_i f(x)|\leq b \leq (\mathbb E [f(G)]^2)^{1/2} \leq C \mathbb E f(G) \quad {\rm a.e.}, \] and the Cauchy-Schwarz inequality we may write
\begin{align*}
\mathbb E |g_i \partial_if(G)| &= \mathbb E |g_i \partial_i f(G)| \mathbf 1_{A^c} + \mathbb E |g_i \partial_i f(G)| \mathbf 1_A \\
&\leq \frac{1}{c_0 \mathbb Ef(G)} \mathbb E [|g_i \partial_i f(G)| f(G)] + b\sqrt{\mathbb P(A)} \\
&\leq \frac{1}{c_0 \mathbb Ef(G)} \frac{\mathbb Ef^2(G)}{n} + \frac{b}{n} \\
&\leq \frac{C_0 \mathbb Ef(G)}{n},
\end{align*} as required.

\end{note}

\section{Concentration for norms with moderate unconditional structure} \label{sec:con-d-unc}

In this section we study Question \ref{que:main} and we prove that it has an affirmative answer for normed spaces
which do not have extremal unconditional basis constant, by establishing Theorem \ref{thm:main-1}. In fact our argument takes into
account a slightly weaker notion, that of the {\it left random unconditional constant} (see Theorem \ref{thm:conc-U} for the details). 
The approach we present uses Proposition \ref{prop:2-level} and the Borsuk-Ulam theorem \cite{Bor} (see also \cite{Mat}). We start with the following:

\begin{lemma} [Balancing the partial derivatives] \label{lem:Borsuk-1} 
Let $f:\mathbb R^m \to \mathbb R$ be a $C^1$-smooth function with bounded partial derivatives and $q>0$. 
Then there exists a diagonal matrix
$\Lambda={\rm diag} (\lambda_1, \ldots, \lambda_m)$ with
\begin{itemize}
\item [\rm (a)] $\|\Lambda\|_{\rm HS} =1$ and
\item [\rm (b)] $\| \partial_i (f\circ \Lambda)\|_{L_q(\gamma_m)}= \|\partial_j (f\circ \Lambda)\|_{L_q(\gamma_m)}$ for $i,j=1,\ldots,m$.
\end{itemize}
\end{lemma}

\noindent {\it Proof.} For each $1\leq j<m$ consider the functions $h_j :S^{m-1}\to \mathbb R$ defined by
\begin{align*}
h_j(\lambda) := \lambda_j \|(\partial_j f)\circ \Lambda\|_{L_q}- \lambda_{j+1} \|(\partial_{j+1} f)\circ \Lambda \|_{L_q},
\end{align*} where $\Lambda= {\rm diag}(\lambda_1, \ldots, \lambda_m)$. The dominated convergence theorem 
and the continuity of $\partial_j f$ imply the continuity of $h_j$, while the symmetry of $\gamma_m$ implies that $h_j$ is odd, that
is $h_j(-\lambda)=-h_j(\lambda)$ for all $\lambda\in S^{m-1}$. 
Hence, if we consider the mapping $H: S^{m-1}\to \mathbb R^{m-1}$ defined by
\begin{align*}
H(\lambda_1, \ldots, \lambda_m) := \big( h_1(\lambda), \ldots, h_{m-1}(\lambda) \big),
\end{align*} we readily see that it is continuous and odd. Therefore, by the Borsuk-Ulam antipodal theorem \cite{Mat} we obtain 
 $\lambda \in S^{m-1}$ such that $H(\lambda)=0$, that is
 \begin{align} \label{eq:6.2}
 \lambda_i \|(\partial_i f) \circ \Lambda\|_{L_q} = \lambda_j \|(\partial_j f) \circ \Lambda\|_{L_q}, \quad i,j=1,2,\ldots,m. 
 \end{align} In particular $\| \partial_i (f\circ \Lambda)\|_{L_q} = \|\partial_j (f\circ \Lambda)\|_{L_q}$ for all $i,j$ which
 proves the assertion.  \prend

 \begin{remarks} \label{rems-borsuk} 1. Note that if $f$ is not constant in any proper subspace, then $\lambda_i>0$ for all $i$. Indeed; note that 
 the set $\sigma :=\{ i : \lambda_i\neq 0 \}$ is non empty. Assuming that $\sigma^c\neq \emptyset$, by \eqref{eq:6.2} we 
 get $\| (\partial_{i} f) \circ \Lambda\|_{L_q} =0 $ for all $i\in \sigma$. Note that $(\partial_{i} f) \circ \Lambda \equiv 0$ for all $i\in \sigma$, by the continuity. 
 It follows that $\partial_{i} f\equiv \bf 0$ on $\Lambda(\mathbb R^m)= \mathbb R^\sigma \equiv [e_i : i\in \sigma]$ for all 
 $i\in\sigma$, i.e. $f|_{\mathbb R^\sigma}$ is constant. Moreover, \eqref{eq:6.2} implies that all $\lambda_j$ have the same sign. Since
 $H(\lambda)=H(-\lambda)=0$ we may assume that $\lambda_j \geq 0$ for all $j$. 
 
 \smallskip
 
\noindent 2. Note that the proof of Lemma \ref{lem:Borsuk-1} can also be applied on the boundary $\partial K$ of any symmetric 
convex body $K$ in $\mathbb R^m$, thus we may also have $\|\lambda\|_K=1$ instead of $\|\lambda\|_2=1$. However, this is not 
crucial for our purposes.
\end{remarks}

\begin{lemma} \label{lem:Borsuk-2}
Let $f: \mathbb R^m\to \mathbb R$ be a $C^1$-smooth Lipschitz map which is not a constant in any proper subspace. 
Then, there exist $\lambda_1,\ldots,\lambda_m > 0$ such that $\sum_{j=1}^m \lambda_j^2=1$ and 
\begin{align*}
\|\partial_j (f\circ \Lambda)\|_{L_1(\gamma_m)} \leq  \frac{1}{m} a(f\circ \Lambda), \quad j=1,2,\ldots, m,
\end{align*} where $\Lambda={\rm diag}(\lambda_1, \ldots, \lambda_m)$. 
\end{lemma} 

\noindent {\it Proof.} Since $f$ is $C^1$-smooth we may consider the diagonal matrix $\Lambda$ 
from Lemma \ref{lem:Borsuk-1} and by taking into Remark \ref{rems-borsuk}.1 we also have $\lambda_j>0$ for all $j$. Note that,
\begin{align*}
\|\partial_j (f\circ \Lambda)\|_{L_1(\gamma_m)} =\frac{1}{m} \sum_{j=1}^m \|\partial_j (f\circ \Lambda)\|_{L_1(\gamma_m)} =
\frac{1}{m} \int_{\mathbb R^m} \| \nabla (f\circ \Lambda )\|_1 \,  d\gamma_m\leq \frac{a(f\circ \Lambda)}{m}, 
\end{align*} for all $j\leq m$, as required. \prend

\medskip

Now we are ready to prove the main result of this section. To this end, let us recall a variant of a one-sided unconditional constant.
The {\it random unconditional divergence constant} of a normed space $X=(\mathbb R^m, \|\cdot\|)$, denoted by ${\rm rud}(X)$, is the least $L>0$ for
which there exists a basis $(x_i)_{i=1}^m$ in $X$ such that 
\begin{align}
\left\| \sum_{i=1}^m \alpha_i x_i \right\| \leq L \mathbb E_\varepsilon \left\| \sum_{i=1}^m \varepsilon_i \alpha_i x_i \right\|,
\end{align} for all scalars $(\alpha_i)_{i=1}^m$. Note that ${\rm rud}(X) \leq {\rm unc}(X)$ (see \cite[Section 6]{FJ} and \cite{AT} for further details). 
With this terminology we have the following.

\begin{theorem} \label{thm:conc-U}
Let $X=(\mathbb R^m, \|\cdot\|)$ be a normed space and let $L={\rm rud}(X)$. Then, there exists $T\in GL(m)$ such that
\begin{align*}
\mathbb P \left(  \big| \|TG\| - \mathbb E\|TG\| \big| > \varepsilon \mathbb E\|TG\| \right) \leq 
C \exp \left(-c\max\{\varepsilon,\varepsilon^2\} \log \left(e+ \frac{m}{L^2} \right) \right),
\end{align*} for all $\varepsilon>0$, where $G$ is the standard $m$-dimensional Gaussian vector and $c,C>0$ are universal constants.
\end{theorem}

\noindent {\it Proof.} After applying an invertible linear transformation we may assume that
\begin{align*} 
\| y\|\leq L \mathbb E_\varepsilon \left\| \sum_{i=1}^m \varepsilon_i y_i e_i \right\|, \quad y\in \mathbb R^m,
\end{align*} where $(e_i)_{i\leq m}$ is the standard basis on $\mathbb R^m$. Equivalently, we have
\begin{align}\label{eq:d-cube}
\sup_{\varepsilon_i=\pm 1} \left\| \sum_{i=1}^m \varepsilon_i y_i e_i\right\| \leq 
L \mathbb E_\varepsilon \left\| \sum_{i=1}^m \varepsilon_i y_i e_i \right\|, \quad y\in \mathbb R^m.
\end{align}
First we consider the case that the given norm $\|\cdot \|$ is smooth. Let 
$\Lambda={\rm diag}(\lambda_1, \ldots, \lambda_m)$ be the diagonal matrix from Lemma \ref{lem:Borsuk-2}. We set 
\begin{align*}
a_\Lambda := \|\Lambda:\ell_\infty^m \to X\|, \quad b_\Lambda:= \|\Lambda: \ell_2^m \to X\|, \quad k_\Lambda:= \frac{(\mathbb E\| \Lambda G\|)^2}{b_\Lambda^2}.
\end{align*} 
Then, by the distributional inequality \eqref{eq:exp-conc} in conjunction with Lemma \ref{lem:Borsuk-2} we obtain
\begin{align} \label{eq:3.4}
\mathbb P\left( \big| \|\Lambda G\| - \mathbb E\|\Lambda G\| \big| > \varepsilon \mathbb E\|\Lambda G\| \right) \leq 
4 \exp \left( -c \varepsilon \sqrt{ k_\Lambda \log \left( e+\frac{m b_\Lambda^2}{a_\Lambda ^2} \right)} \right), \quad \varepsilon >0.
\end{align}  Employing \eqref{eq:d-cube} we arrive at the following estimate:
\begin{align*}
a_\Lambda = \max_{\varepsilon_i = \pm 1} \left\| \sum_{i=1}^m \varepsilon_i \lambda_i e_i \right\| \leq 
L \mathbb E_\varepsilon \left\| \sum_{i=1}^m \varepsilon_i \lambda_i e_i \right\| \leq L\sqrt{\frac{\pi}{2}} \mathbb E\|\Lambda G\|,
\end{align*} where in the last passage we have used the contraction principle (see e.g. \cite[Chapter 4]{LedTal-book} or \cite[Proposition 1]{Pis-type}).
Plugging the above estimate into \eqref{eq:3.4} we obtain
\begin{align*}
\mathbb P\left( \big| \|\Lambda G\| - \mathbb E\|\Lambda G\| \big| > \varepsilon \mathbb E\|\Lambda G\| \right) \leq 
4 \exp \left(-c \varepsilon \sqrt{k_\Lambda \log \left( e+ \frac{ m}{ L^2 k_\Lambda}  \right) } \right),
\quad \varepsilon >0.
\end{align*} Applying the latter for $\varepsilon \simeq 1$ and taking into account \eqref{eq:ldr-lb} we readily see that 
$k_\Lambda \geq c \log(e+m/L^2)$. This proves the result in the smooth case.

For the general case, recall that for the given norm $\|\cdot\|$ and for any $\delta\in (0,1)$ there exists a smooth norm $\|\cdot\|_\delta$
such that
\begin{align} \label{eq:smooth}
(1-\delta) \|x\| \leq \|x\|_\delta \leq (1+\delta) \|x\|, 
\end{align} for all $x\in \mathbb R^m$, see e.g. \cite{Schn-book}. We fix $0<\delta \leq (7+\log m)^{-1} $ 
and we apply the result for $\|\cdot\|_\delta$,
thus we get $T=T_\delta \in GL(m)$ such that
\begin{align*}
\mathbb P \left( \big| \|TG\|_\delta - \mathbb E\|TG\|_\delta \big| > \varepsilon \mathbb E\|TG\|_\delta \right)
\leq 4e^{-c \varepsilon \log (e+ m/ L_\delta^2) }, 
\end{align*} for all $\varepsilon > 0$, where $L_\delta={\rm rud}(X_\delta)$ for $X_\delta=(\mathbb R^m, \|\cdot\|_\delta)$.
One may check, using \eqref{eq:smooth}, that for $\varepsilon> 8\delta $ we have
\begin{align*}
\mathbb P \left(  \big| \|TG\| - \mathbb E\|TG\| \big| > \varepsilon \mathbb E\|TG\| \right) \leq 
\mathbb P \left(  \big| \|TG\|_\delta - \mathbb E\|TG\|_\delta \big| > \frac{\varepsilon}{2} \mathbb E\|TG\|_\delta \right). 
\end{align*} Thus, we obtain
\begin{align*}
\mathbb P \left(  \big| \|TG\| - \mathbb E\|TG\| \big| > \varepsilon \mathbb E\|TG\| \right) \leq 4 e^{- \frac{c}{2} \varepsilon \log(e+m/L_\delta^2)} ,
\end{align*} for all $\varepsilon > 8\delta$. By adjusting the universal constants the previous estimate holds true for all $\varepsilon>0$. 
On the other hand we may easily check that 
$$L_\delta \leq \frac{1+\delta}{1-\delta} L \leq 2 L.$$
The proof is complete.  \prend

\begin{remarks} 1. It is somewhat unexpected that the almost optimal concentration is established with an isomorphic parameter.
This suggests that Question \ref{que:main} seems plausible to have an affirmative answer and the above partial result consists of a deficiency
of the approach, which uses the $L_1-L_2$ estimate for the quantification of the problem. In fact, in the next paragraph we show how one can
eliminate ${\rm rud}(X)$ and overcome this obstacle.

\smallskip

\noindent 2. The {\it random unconditional convergence constant} ${\rm ruc}(X)$ is defined similarly as the least $R>0$ for which there
exists a basis $(x_i)$ of $X=(\mathbb R^n, \|\cdot\|)$ such that
\begin{align*}
\mathbb E \left \|\sum_{i=1}^n \varepsilon_i \alpha_i x_i \right\| \leq R \left \| \sum_{i=1}^n \alpha_i x_i \right\|,
\end{align*} for all scalars $(\alpha_i)_{i\leq n}$. While RUC bases have been previously studied (see e.g. \cite{BKPS}, \cite{AT}) and their 
extremal asymptotic behavior was established in \cite{Ba-wrgl}, it was not until recently that its left analogue, i.e. ${\rm rud}(X)$, was 
put forward in systematic study; see e.g. \cite{AT}. Consulting specialists
we couldn't locate a precise reference for the study of the latter in the context of high-dimensional normed spaces. 
These notions and more will be part of a detailed study which will appear elsewhere.

\end{remarks}

\section{Probabilistic dichotomy and Dvoretzky's theorem} \label{sec:dicho}

In this Section we prove Theorem \ref{thm:main-2} and its corollaries mentioned in the Introduction.
Our first main ingredient is the classical Dvoretzky-Rogers lemma from \cite{DR}. Another crucial tool in our approach is the Alon-Milman theorem. 
The idea to use this dichotomy in this problem can be traced back to the work of Schechtman \cite{Sch-eps-exist}.

\begin{lemma} [Dvoretzky-Rogers] \label{lem:DR}
Let $X=(\mathbb R^n, \|\cdot\|)$ be a normed space for which $B_X$ is in John's position. Then, there exists 
an orthonormal basis $v_1,\ldots, v_n$ such that 
\begin{align*}
1=\|v_k\|_2\geq \|v_k\|\geq \sqrt{1-\frac{k-1}{n} },
\end{align*} for $k=1,2,\ldots, n$. In particular, $\|v_j\|\geq 1/\sqrt{2}$ for $j=1,\ldots, \lfloor n/2 \rfloor$.
\end{lemma}

\begin{remark} [W.B. Johnson] \label{rem:Johnson} 
Starting with the above orthonormal basis one may redefine the vectors to get a new 
orthonormal basis $(w_i)$ with $\|w_i\| \geq 1/4$ for all $i \leq n$. This remark is due to Bill Johnson \cite{Bill}. 
We would like to thank him for allowing us to include his elegant argument here. Assume for simplicity that 
$n=2s$. For each $i=1,\ldots, s$, if $\|u_{s+i}\|\geq 1/4$ we set $w_{s-i+1}=u_{s-i+1}$ and 
$w_{s+i}=u_{s+i}$, while if $\|u_{s+i}\|<1/4$ we replace $w_{s-i+1}=\frac{u_{s-i+1}+u_{s+i}}{\sqrt{2}}$ and $w_{s+i}=\frac{u_{s-i+1} -u_{s+i} }{\sqrt{2}}$. 
Note that the vectors $(w_i)_{i\leq n}$ are still orthonormal and $\|w_{s-i+1}\| , \|w_{s+i}\| \geq 1/4$. Indeed; by construction we have
\[ \min\{ \|w_{s-i+1}\|, \|w_{s+i}\| \} 
\geq \frac{1}{\sqrt{2}} ( \|u_{s-i+1}\| - \|u_{s+i}\|) > \frac{1}{\sqrt{2}} \left( \frac{1}{\sqrt{2}} -\frac{1}{4} \right) >\frac{1}{4} , \]
where we have used that $\|u_{s-i+1}\| \geq 1/\sqrt{2}$ for all $i\leq s$ by Lemma \ref{lem:DR}.  \prend
\end{remark}

We will also need the following theorem of Alon and Milman \cite{AM}; see also \cite{Tal-cube} for an alternative
simpler proof.

\begin{theorem}[Alon-Milman, Talagrand] \label{thm:AM}
Let $X$ be a normed space and let $T:\ell_\infty^n\to X$. We set 
\[ a=\|T : \ell_\infty^n \to X \| \quad and \quad M_n= \mathbb E_\varepsilon \left\| \sum_{i=1}^n \varepsilon_i Te_i \right\| .\]
Assuming that $\|Te_i\|\geq 1$ for all $i$, there exists $\sigma\subset [n]$ with 
$|\sigma|\geq cn/a$ such that
\begin{align*}
\frac{1}{2} \max_{i\in \sigma} |\alpha_i| \leq \left\| \sum_{i\in \sigma} \alpha_i Te_i \right\| \leq 4 M_n \max_{i\in \sigma}|\alpha_i|,
\end{align*} for all scalars $(\alpha_i) \subset \mathbb R$.
\end{theorem}

Alon-Milman's proof yields $\sigma \subset [n]$ with $|\sigma| \geq cn^{1/2}/ M_n$. The improved 
estimate stated above is due to Talagrand. We are now ready to prove the key result of this section.
Note that Theorem \ref{thm:main-2} will follow from the next result.

\begin{theorem} \label{thm:dicho-1}
Let $X=(\mathbb R^n, \|\cdot\|)$ be a normed space for which $B_X$ is in John's position and let $0< \delta <1/2$. 
Then, at least one of the following conditions holds:
\begin{itemize}
\item Either $k(X) \geq n^{1/2-\delta}$,

\item Or there exists a subspace $F$ with $\dim F=m\geq c n^{1/2}$ and a
linear isomorphism $T:F \to F$ such that for all $t>0$,
\begin{align*}
\mathbb P \left( \big| \|TZ\| - \mathbb E\|TZ\| \big| >t\mathbb E\|TZ\| \right) \leq C e^{-c\delta \max \{ t^2,t \} \log m}, \quad Z\sim N({\bf 0},I_F),
\end{align*} where $c,C>0$ are universal constants.
\end{itemize}
\end{theorem}

\noindent {\it Proof.} Fix $0<\delta<1/2$. Let $(w_j)$ be an orthonormal basis with $\|w_j\|\geq 1/4$ for all $j$ and let
$k(X) \leq n^{1/2-\delta}$. Then, we may write
\begin{align*}
\sqrt{ \frac{2}{\pi} } \mathbb E\left\| \sum_{i=1}^n \varepsilon_i w_i \right\| \leq \mathbb E\left\| \sum_{i=1}^n g_i w_i \right\|
= \mathbb E\|G\|= \sqrt{k(X)},
\end{align*} where in the first inequality we have used the contraction principle, see \cite[Chapter 4]{LedTal-book} or \cite[Proposition 1]{Pis-type}.
Using Theorem \ref{thm:AM} we obtain a subset 
$\sigma \subset [n]$ with $|\sigma| \geq c n /a\geq c\sqrt{n}$, where $a=\|i: \ell_\infty^n \to X\|$ and
\begin{align*}
\frac{1}{8} \max_{i\in \sigma} |\alpha_i| \leq \left\| \sum_{i\in \sigma} \alpha_i w_i \right\| \leq  4 M_n \max_{i\in \sigma} |\alpha_i| ,
\end{align*} for all $(\alpha_i)_{i \in \sigma} \subset \mathbb R$, where $M_n=\mathbb E\|\sum_i \varepsilon_iw_i\|$. 
Note that the subpace $(F, \|\cdot\|)$ with $F={\rm span} \{w_i : i\in \sigma\}$ 
satisfies $d(F, \ell_\infty^\sigma)\leq 32 M_n \leq C \sqrt{k(X)}$. 
Thus, by Theorem \ref{thm:conc-U} there exists a linear isomorphism $T:F\to F$ such that 
\begin{align*}
\mathbb P \left( \big| \|TZ\| -\mathbb E \|TZ\| \big| >t \mathbb E \|T Z\| \right) &\leq 
Ce^{-c \max\{ t^2 , t \} \log(e +c|\sigma|/ k(X)) } \\
&\leq Ce^{-c' \delta \max\{ t^2, t \} \log |\sigma|},
\end{align*} for all $t >0$. The proof is complete. \prend

\medskip

In the following theorem, which is immediate consequence of Theorem \ref{thm:dicho-1}, we summarize Theorem \ref{thm:main-3} and 
Corollary \ref{cor:main-4}. We would like to thank K. Tikhomirov who kindly pointed out to us that the linear map in Theorem \ref{thm:dicho-1} can
be lifted up to an invertible linear map, whence obtaining an affirmative answer to Question \ref{que:main} as it is stated.

\begin{theorem} 
Let $X=(\mathbb R^n, \|\cdot\|)$ be a normed space. Then, there exists a linear transformation $S\in GL(n)$
such that for all $t>0$, 
\begin{align*}
\mathbb P\left( \big| \|SG\| - \mathbb E\|SG\| \big| > t \mathbb E\|SG\| \right) \leq 
Ce^{-c\max\{t^2, t \} \log n}, \quad G\sim N({\bf 0},I_n).
\end{align*} In particular, for any $\varepsilon \in (0,1)$ the random $k$-dimensional section of $S^{-1}B_X$ with 
$k \leq c\varepsilon \log n /\log(1/\varepsilon)$ is $(1+\varepsilon)$-spherical with probability greater than $1-n^{-c\varepsilon}$.
\end{theorem}

\noindent {\it Proof.} We may assume that $B_X$ is in John's position and $k(X)\leq n^{1/3}$, otherwise there is nothing to prove. 
Then, Theorem \ref{thm:dicho-1} yields the existence of a subspace $F$ with $\dim F=m \simeq \sqrt{n}$ and $T : F\to F$ with $T \in GL(F)$ such that for every $\varepsilon > 0$,
\begin{align*}
\mathbb P \left( \big| \|TZ\| -\mathbb E\|TZ\| \big| > \varepsilon \mathbb E\|TZ\| \right) \leq C \exp \left( -c\varepsilon \log n \right),  
\quad Z\sim N({\bf 0},I_F).
\end{align*} 
Let $S : \mathbb R^n \to \mathbb R^n$ be the operator defined by
\begin{align*}
S(x,y) = Tx + \lambda y, \quad x\in F, \; y\in F^\perp, \quad \lambda= \frac{\mathbb E\|TZ\| }{\mathbb E\|W\|_2 \log n}, \quad W\sim N({\bf 0},I_{F^\perp}).
\end{align*} Note that by the contraction principle, the triangle inequality and the fact that 
$\|\cdot\| \leq \|\cdot\|_2$ we have
\begin{align} \label{eq:exp-lift}
\mathbb E\|TZ\| \leq \mathbb E\|SG\|  \leq \left( 1+ \frac{1}{ \log n}\right) \mathbb E\|TZ\|, \quad G\sim N({\bf 0},I_n), \; Z\sim N({\bf 0},I_F).
\end{align} Next, we have the following:

\smallskip

\noindent {\it Claim.} For all $\varepsilon >0 $ we have
\begin{align} \label{eq:global-con}
\mathbb P \left( \big| \|SG\| -\mathbb E \|SG\| \big | > \varepsilon \mathbb E\|SG\| \right) \leq Ce^{-c\varepsilon \log n}.
\end{align}

\smallskip

\noindent {\it Proof of Claim.} Let $Z\sim N({\bf 0},I_F)$, $W\sim N({\bf 0},I_{F^\perp})$ and $G=Z+W$. Then, we may write
\begin{align*}
\mathbb P\left( \|SG\| > (1+\varepsilon) \mathbb E\|SG\| \right) &\leq \mathbb P\left( \|TZ\| > (1+\varepsilon) \mathbb E\|SG\| - \lambda \|W\| \right) \\
&\leq \mathbb P \left( \|TZ\| >(1+\varepsilon) \mathbb E\|TZ\| - 10 \lambda \mathbb E\|W\|_2 \right) + \mathbb P \left(  \|W\| > 10 \mathbb E\|W\|_2 \right) \\
&\leq \mathbb P \left( \|TZ\| > (1+\varepsilon)\mathbb E \|TZ\|- \frac{10 \mathbb E\|TZ\|}{ \log n} \right) + \mathbb P \left( \|W\|_2 > 10 \mathbb E\|W\|_2 \right) \\
&\leq \mathbb P \left( \|TZ\| > \left (1+\frac{\varepsilon}{2}\right)\mathbb E \|TZ\| \right) +  \mathbb P \left( \|W\|_2 > 10 \mathbb E\|W\|_2 \right) \\
&\leq Ce^{-c\varepsilon \log n},
\end{align*} for all $\frac{20}{\log n}<\varepsilon <1$, where we have also used the fact that 
\begin{align*}
\mathbb P (\|W\|_2 >10 \mathbb E\|W\|_2 ) \leq  e^{-cn} , \quad W \sim N({\bf 0},I_{F^\perp}).
\end{align*} For the deviation below the mean we may write
\begin{align*}
\mathbb P( \|SG\| < (1-\varepsilon) \mathbb E\|SG\| ) &\leq \mathbb P \left(  \|TG\| < (1-\varepsilon)\left(1+\frac{1}{\log n} \right) \mathbb E\|TG\| \right) \\
&\leq \mathbb P \left( \|TG\| < \left( 1-\frac{\varepsilon}{2} \right) \mathbb E\|TG\| \right) \\
&\leq Ce^{-c\varepsilon \log n},
\end{align*} for all $\frac{2}{\log n}<\varepsilon<1$, where we have used \eqref{eq:exp-lift}. 
The claim follows for all $\varepsilon \in (0,1)$ by adjusting the universal constants. Note that the 
estimate \eqref{eq:global-con}, combined with \eqref{eq:ldr-lb}, yields that $k(S^{-1}(B_X)) \geq c\log n$, hence the two-level tail estimate readily follows. \prend

\section{Further remarks and questions} \label{sec:rems}

We end this note with some concluding comments that arise from our work.

\medskip

\noindent {\bf \S \, 1. Concentration for norms close to $\ell_2$.} Here we recall the fact that norms close to the $\ell_2$ norm share analogous concentration, as was 
claimed in the Introduction. Let $\alpha\geq 1$ and let $\|\cdot \|$ be a norm on $\mathbb R^n$ such that \[ \|x\|_2\leq \|x\|\leq \alpha \|x\|_2, \quad x\in \mathbb R^n. \] 
Note that $x\mapsto \|x\|$ is $\alpha$-Lipschitz and $\mathbb E\|G\| \geq \mathbb E\|G\|_2 \simeq \sqrt {n}$, hence the standard Gaussian 
concentration inequality \eqref{eq:con-Gauss} implies
\[ \mathbb P \left( \big| \|G\| - \mathbb E\|G\| \big| >\varepsilon \mathbb E\|G\| \right) \leq C\exp(-c \varepsilon^2 n/ \alpha^2), \quad \varepsilon>0. \]
Moreover, one can show that this estimate is essentially optimal, i.e.
\[\mathbb P \left( \big| \|G\| - \mathbb E\|G\| \big| >\varepsilon \mathbb E\|G\| \right) \geq c\exp(-C \varepsilon^2 n \alpha^2), \quad \varepsilon>0. \] For a proof of the 
latter estimate the reader is referred to \cite[Lemma 6.1]{PVZ}.

\medskip

\noindent {\bf \S \, 2. Concentration in John's position.} We provide an explicit construction of norms 
which shows that in John's position the concentration estimate \eqref{eq:con-norms} cannot be improved.

\begin{proposition}
Let $X=(\mathbb R^n, \|\cdot\|)$ be a normed space and let $Y:= (X \oplus \ell_2^m)_\infty$, i.e. 
\begin{align*}
\|y\|_Y= \max \left\{ \|x\| , \|z\|_2\right\}, \quad y=(x; z) \in \mathbb R^n\times \mathbb R^m. 
\end{align*} Suppose that $B_X$ is in John's position. Then, we have the following:
\begin{itemize}
\item [\rm i.] $B_Y$ is also in John's position and for $m\geq Ck(X)$ 
we have ${\rm Var}\|G\|_Y \geq c b(Y)^2$.
\item [\rm ii.] For $m\simeq k(X)$, the norm $\|\cdot\|_Y$ exhibits the following concentration
\begin{align*}
ce^{-Ct^2k(X)} \leq \mathbb P \left( \big| \|G\|_Y - \mathbb E\|G\|_Y \big|>t \mathbb E\|G\|_Y \right) \leq Ce^{-ct^2 k(X)}, \quad t>0.
\end{align*}
\item [\rm iii.] For any $\varepsilon\in (0,1)$ we have $k_r(Y,\varepsilon) \simeq \varepsilon^2 k(X)$.
\end{itemize}
\end{proposition}

\noindent {\it Proof.} i. It is easy to verify that $B_Y$ is in John's position. Indeed; for all $y=(x;z)$ we clearly have
\begin{align*}
\|y\|_Y = \|(x;z)\|_Y \leq \max\{\|x\|_2,\|z\|_2 \}\leq \sqrt{\|x\|_2^2+ \|z\|_2^2}=\|(x; z)\|_2.
\end{align*} Furthermore, there exist $u_1,\ldots,u_s\in S^{n-1}$ contact points, i.e. $\|u_i\|=\|u_i\|_\ast=1$ 
and $c_1,\ldots c_s>0$ such that $I_{\mathbb R^n}=\sum_{j=1}^s c_j u_j\otimes u_j$.
Hence, we have
\begin{align*}
I_{\mathbb R^n\times \mathbb R^m} = 
\sum_{i=1}^s c_i (u_i;\mathbf 0_{\mathbb R^m})\otimes (u_i;\mathbf 0_{\mathbb R^m}) +\sum_{j=1}^m e_{n+j}\otimes e_{n+j}.
\end{align*} By the converse of John's theorem (see \cite{Ba} for a proof) we conclude that $B_Y$ is in John's position.
Let $Z,W$ be independent Gaussian vectors with $Z\sim N({\bf 0},I_m)$ and $W\sim N({\bf 0},I_n)$ and let $G=(W,Z)\sim N({\bf 0}, I_{n+m})$. 
If $A=\{y=(x; z)\in \mathbb R^{n+m} : \|x\| \leq \|z\|_2\}$, then we may write
\begin{align*}
{\rm Var}[ \|G\|_Y] & \geq \frac{1}{2} \iint_{A\times A} ( \|z\|_2-\|z'\|_2 )^2 \, d\gamma_{n+m}(y) d\gamma_{n+m}(y') \\
& \geq {\rm Var}[\|Z\|_2 ] - \frac{1}{2} \iint_{(A\times A)^c} (\|z\|_2 -\|z'\|_2)^2 \, d\gamma_{n+m}(y) \, d\gamma_{n+m}(y') \\
&\geq c_0 - C_0 \sqrt{ \mathbb P \left( (A\times A)^c \right) },
\end{align*} where we have used the Cauchy-Schwarz inequality and the fact that
\begin{align*}
c_0 \leq {\rm Var}(\|Z\|_2) \leq \left( \mathbb E \big| \|Z\|_2 - \|Z'\|_2 \big|^4 \right)^{1/2} \leq C_0,
\end{align*} where $Z'$ is an independent copy of $Z$ (see e.g. \cite[Proposition 4.4]{PVZ}). On the other hand we have
\begin{align*}
\mathbb P \left(  (A\times A)^c \right) &\leq 2 \mathbb P( \|Z\|_2 < \|W\|) \\
&\leq 2 \left[ \mathbb P \left( \|Z\|_2 \leq \delta \sqrt{m} \right) + \mathbb P \left( \|W\| > \delta \sqrt{m} \right) \right] \\
&\leq (c_1\delta)^m + c_2 e^{-c_3\delta^2 m},
\end{align*} provided that $\delta\sqrt{m}\geq 2\mathbb E\|W\| = 2\sqrt{k(X)}$. Choosing $\delta$ sufficiently small universal constant 
we obtain $\mathbb P((A\times A)^c)\leq e^{-c m}$, hence for $m\geq C \delta^{-2} k(X)$ we get
${\rm Var} (\|G\|_Y) \geq c_0'$, as required. 

\smallskip

\noindent ii. Note that $k(Y)\simeq \max\{ k(X), m \}$. Taking 
into account (i) and employing the main result of \cite{V-tight} we get the assertion.

\noindent iii. Recall that $k_r(Y,\varepsilon)$ is the maximal $k$ for which the random $k$-dimensional subspace of $Y$ 
is $(1+\varepsilon)$-Euclidean with probability at least 2/3. Once we have established the sharp concentration for the norm, 
it is routine to check that $k_r(Y, \varepsilon)\simeq \varepsilon^2 k(X)$. For the details see \cite[Section 3]{V-tight}. \prend

\begin{remarks} 1. The above construction, when $X=\ell_\infty^n$, yields a 1-unconditional norm for which 
the concentration estimate \eqref{eq:con-John} in John's position cannot be improved. Essentially this example is 
due to K. Tikhomirov, who proves part (iii) in \cite{Tik-unc} since his focus is on the dependence on $\varepsilon$ 
in Dvoretzky's theorem. His approach is completely different from the one presented here and lies in delicate estimates 
for singular values of Gaussian matrices. 

\smallskip

\noindent 2. Applying the above construction for $X=\ell_q^n, \; 2\leq q \leq \infty$ and $m\simeq k(\ell_q^n)$, we get spaces 
which are in John's position, exhibit optimal concentration in terms of the Lipschitz constant and have Dvoretzky number of 
all possible range, i.e. $k(Y)\simeq k(\ell_q^n) \simeq k(\ell_q^N) \in (\log N, N)$.

\smallskip

\noindent 3. The same spaces $Y$ as above yield examples of 1-unconditional normed spaces which are in John's position, 
are of cotype $q$ with constant $C_q(Y) \simeq C_q(\ell_q^n)$ and satisfy optimal concentration as in \eqref{eq:con-norms}. 
This shows that the consideration of \cite{FLM} on the random version of Dvoretzky's theorem for spaces with cotype in John's position 
cannot be improved. This also shows that their approach to study the corresponding question for subspaces of $L_q, 2<q<\infty$,
by viewing them as spaces with cotype $q$, is insufficient and other tools are required to obtain the optimal estimates, 
see \cite{PV-dvo-Lp} for the details.

\end{remarks}

\noindent {\bf \S \, 3. Hypercontractive measures.} It is worth mentioning that Proposition \ref{prop:2-level} holds true 
for any hypercontractive measure, since such measures satisfy \eqref{eq:Tal-L1-L2}; see \cite{CorLed}. Recall that a measure $\mu$
is said to be hypercontractive with constant $\rho$ if it satisfies a log-Sobolev inequality with constant $\rho>0$, i.e.
\[{\rm Ent}_\mu(f^2)= \mathbb E_\mu f^2\log f^2 - \mathbb E_\mu f^2 \log \mathbb E_\mu f^2 \leq \frac{2}{\rho} \mathbb E_\mu\|\nabla f\|_2^2, \]
for any smooth function $f$. In particular, we have the following: 
Let $\mu$ be a hypercontractive and symmetric Borel probability measure on $\mathbb R^n$ with constant $\rho>0$, 
and let $\|\cdot\|$ be an arbitrary norm on $\mathbb R^n$. 
Then, for any smooth Lipschitz map $f:\mathbb R^n\to \mathbb R$ there exists $\lambda= (\lambda_1, \ldots, \lambda_n) \in \mathbb R^n$ with
$\| \lambda\|=1$ such that
\begin{align*}
\mu \left( z\in \mathbb R^n : \big| f(\Lambda z) -\mathbb E_\mu f(\Lambda z) \big| > t \right) \leq 
4 \exp \left( -\frac{ct}{b_\Lambda} \sqrt{\rho \log \left( e+ \frac{n b_\Lambda^2}{a_\Lambda^2}\right) } \right), \quad  t>0,
\end{align*} where $\Lambda={\rm diag}(\lambda_1, \ldots, \lambda_n)$ and $b_\Lambda=b(f\circ \Lambda)$, $a_\Lambda=a(f\circ \Lambda)$.

In addition, if $f$ is a norm, one may get the following variant of Theorem \ref{thm:conc-U}: For any normed space $X=(\mathbb R^n, \|\cdot \|)$ 
there exists a non-singular matrix $T$ such that
\begin{align*}
\mu \left(z : \left| \|Tz\| -\mathbb E_\mu \|Tz\| \right| >t \mathbb E_\mu\|Tz\| \right) \leq  
C\exp\left( -c t\sqrt{ \rho k \log \left( \frac{en}{d(X,\ell_\infty^n)^2} \right) } \right), \quad t>0,
\end{align*} where $k=(\mathbb E_\mu\|Tz\| / \max_{\|z\|_2\leq 1} \|Tz\|)^2$. To this end, one needs to invoke the following elementary fact.

\begin{fact}
If the norm $\|\cdot\|$ on $\mathbb R^n$ satisfies $\|x\|_\infty \leq \|x\|\leq a\|x\|_\infty$ for all $x$, then for any diagonal matrix
$\Lambda$ we have
\begin{align*}
\max_{\|x\|_\infty \leq 1} \|\Lambda x\| \leq a \max_{ \|x\|_2 \leq 1} \|\Lambda x\| .
\end{align*}
\end{fact}

\noindent The details are left to the interested reader.

\medskip

\noindent {\bf \S \, 4. On the parameter $\beta$.} The following parameter, referred to as the {\it normalized variance} is 
introduced in \cite{PV-sdi} (see also \cite{PV-var}) for the study of sharp Gaussian small deviation inequalities and small ball probabilities
for norms. For any normed
space $X=(\mathbb R^n,\|\cdot\|)$ we define
\begin{align*}
\beta(X) =\beta(B_X) =\frac{{\rm Var}(\|G\|)}{(\mathbb E\|G\|)^2}, \quad G\sim N({\bf 0},I_n).
\end{align*} It is also known (see e.g. \cite{PV-var}) that $\beta(X) \geq \beta(\ell_2^n)\simeq 1/n$. We define further
\begin{align*}
\mathfrak B(X) :=\min_{T\in GL(n)} \beta(TB_X).
\end{align*} 
In \cite{PV-dvo-Lp} we prove that for any $n$-dimensional subspace $X$ of $L_p, \; 1\leq p<\infty$ one has
\begin{align*}
\mathfrak B(X) \leq \frac{e^{cp}}{n},
\end{align*} which is clearly of minimal possible order (up to constants of $p$). In \cite{PVZ} and \cite{LT} the parameter $\beta(\ell_p^n)$ is
estimated asymptotically with respect to $n$ and $p$ (when $p$ grows along with $n$).

In the light of Proposition \ref{prop:con-unc} we get for any 1-unconditional normed 
space $X$ in position of minimal $M$, that
\begin{align*}
\beta (X) \leq \frac{C}{k\log (en/k)}, \quad k=k(X) \geq c\log n,
\end{align*} which is clearly optimal for $X=\ell_\infty^n$. In particular, 
\begin{align*}
\mathfrak B (X) \leq \frac{C}{ (\log n)^2},
\end{align*} for any 1-unconditional normed space $X$. The latter is also derived by K. Tikhomirov in \cite{Tik-unc}.

The main result of Section \ref{sec:dicho} shows moreover, that for any normed space $X=(\mathbb R^n, \|\cdot \|)$ one has 
\begin{align*}
\mathfrak B(X) \leq \frac{C }{(\log n)^2 }.
\end{align*}

\bigskip

\noindent {\bf Acknowledgments.} The authors are grateful to Bill Johnson, Emanuel Milman, Gideon Schechtman and Konstantin Tikhomirov for useful comments.
Thanks also go to the anonymous referees whose valuable suggestions improved the presentation of the paper.
Part of this work was carried out when the second named author was visiting Texas A\&M University on the occasion of the Workshop in Analysis 
and Probability, July 2017. He would like to thank the organizers for the hospitality and the excellent research environment provided.
This material is also based upon work supported by the National Science Foundation under Grant DMS-1440140 while the authors
were in residence at the Mathematical Sciences Research Institute Berkeley, California, during the Fall 2017 semester, on the occasion 
of the Program Geometric Functional Analysis and Applications. The authors would like to thank the faculty of MSRI and the organizers of the program for the exceptional research conditions.

{\footnotesize

\bibliographystyle{alpha}
\bibliography{dichotomy-ref}
}

\medskip

\vspace{.5cm} \noindent 

\begin{minipage}[l]{\linewidth}
  Grigoris Paouris: {\tt grigoris@math.tamu.edu}\\
  Department of Mathematics, Mailstop 3368\\
  Texas A \& M University\\
 College Station, TX 77843-3368\\
  
  \medskip
  
  Petros Valettas: {\tt valettasp@missouri.edu}\\
  Mathematics Department\\
  University of Missouri\\ 
  Columbia, MO 65211\\

\end{minipage}

\end{document}